\documentclass[ 11pt] {article}
\usepackage{a4,  psfrag, epsfig, epsf, amscd, amsfonts, amsmath, amssymb, amstext, amsthm, enumerate}
\begin{document}

\newtheorem{satz}{Theorem}[section]
\newtheorem{defin}[satz]{Definition}
\newtheorem{bem}[satz]{Remark}
\newtheorem{hl}[satz]{Proposition}
\newtheorem{lem}[satz]{Lemma}
\newtheorem{ko}[satz]{Corollary}

\title{If all geodesics are closed on the projective plane }
\author{Christian Pries}
\maketitle

\begin{abstract}
Given a $C^{\infty}$  Riemannian metric $g$ on $\mathbb{R}P^2$
we prove that  $( \mathbb{R}P^2, g)$ has constant curvature 
  iff  all geodesics are closed. Therefore $ \mathbb{R}P^2$ is the first non trivial  example of a manifold such that the smooth Riemannian metrics which involve that all geodesics are closed are unique up to isometries and scaling. This remarkable phenomenon is not true on the 2-sphere, since there is a large set of $C^{\infty}$ metrics whose geodesics are all closed and have the same period 2$\pi$ (called Zoll metrics), but no metric of this set can be obtained from another metric of this set via an isometry and scaling. As a corollary we conclude that all two dimensional P-manifolds
are SC-manifolds.
\end{abstract}

\setcounter{section}{-1}
\section{Introduction}
  Gromoll and Grove proved in [G]
that if $(S^2, g)$ is a P-manifold then it is a SC-manifold (see defintion below).
We prove that this even holds if 
$\mathbb{R}P^2$ is a P-manifold. Moreover we show that $\mathbb{R}P^2$ has constant curvature iff all geodesics are closed. The first result in this direction was proven by Green in [Gr]. He proved that $S^2$ has constant curvature iff it is a Blaschke manifold.
From this theorem it follows that  $(\mathbb{R}P^2, g)$ has constant curvature iff all geodesics are closed, have the same period and are without selfintersections, since the orientable double cover of $(\mathbb{R}P^2, g)$  is then a Blaschke manifold. 
In the complete paper all geodesics are parametrized by arc length and the geodesic flow is complete. The manifolds and the  Riemannian metrics  are $C^{\infty}$.
$\pi: TM\to M$ denotes the canonical projection. $\gamma_v$ denotes the geodesic with
$ \dot{\gamma}_v(0)= v$.
For interested readers who want to know more about P-manifolds
 we refer to the book  [Bs].

\section[]{Preliminaries}

\begin{defin}
$(M,g)$ is called a P-manifold if all geodesics are closed. $(M,g)$ is called a C-manifold if all geodesics are closed and they  have all the same period. $(M,g)$ is called a SC-manifold if $(M,g)$ is  a C-manifold and  all
geodesics are simple closed. The associated metrics are called P-metrics (resp. C-metrics and SC-metrics).
\end{defin}

\begin{defin}
Let $(X,d)$ be a metric space. A dynamical system $\Phi : \mathbb{R} \times X \to X$ is called  
\begin{enumerate}[i)]
\item equicontinuous  (regular) if for all $\epsilon > 0$ there exists a $\delta (\epsilon) > 0$ such that for all $x,y$ with 
$d(x,y) < \delta (\epsilon) $ we have $d(\Phi_t (x), \Phi_t (y))< \epsilon$ for all $t \in \mathbb{R}$.
\item uniformly almost periodic if there exists
for every $\epsilon > 0$ a  $\tau >0 $ such that in every interval $I$ of length $\tau$
there exists a $t\in I$ such that for all $x $ we have $d(\Phi_t (x), x) < \epsilon $. 
\item periodic if there exists an $L > 0$ such that $\Phi_L = Id$, i.e. every point is periodic
and all points have a common period.
\end{enumerate}
\end{defin}

Note that the terms above are independent from the metric which defines the topology if the space $X$ is compact.

\begin{satz}
An equicontinuous flow $\Phi : X \times \mathbb{R}   \to X$ on  a compact metric  space  is uniformly almost periodic.
\end{satz}
Proof: See theorem 2.2 in [A]. {\hfill $\Box$ \vspace{2mm}}

\begin{satz}[Wadsley]
Given a P-manifold $(M,g)$ and let  denote $\Phi:\mathbb{R} \times SM \to SM$ its geodesic flow. Define for $v\in SM$ $P(v)$ to be the least period of $\Phi_t(v)$.
Then $M$ is compact, the geodesic flow is periodic and $P$ is a  discrete, bounded function. Moreover, $P(v)$ is a multiple of the smallest period, i.e. $P(SM)\subset \mathbb{N} P(v_0)$ for some $v_0\in SM$.
\end{satz}
Proof: The compactness of $M$ is clear. Wadsley showed in [W] a characterization of periodic flows on manifolds. One can conclude from this article that the geodesic flow is periodic. A proof of this fact can be found in [Bs] (see lemma 7.11). \\ Since the geodesic flow is periodic (M is a P-manifold), all geodesics have a common period and $P$ is a discrete, bounded function. It is known that in this case 
 $P(v)$ is a multiple of the smallest period, i.e. $P(SM)\subset \mathbb{N} P(v_0)$ for some $v_0\in SM$ (see page 490 in [Be]).{\hfill $\Box$ \vspace{2mm}}

\begin{lem}
If $(M,g)$ is a P-manifold then the geodesic flow $\Phi:\mathbb{R} \times SM \to SM$ 
is equicontinuous with respect to the Sasaki-metric.
\end{lem}

Proof: We know from theorem [1.4] that the flow is periodic and $M$ is compact. Let $L$ denote the  period of the geodesic flow. Choose for $\epsilon > 0$
a $\delta > 0$ such that $d(v,w)< \delta $ for $v,w\in SM$ implies 
$$d(\Phi_t(v),\Phi_t(w) ) < \epsilon $$
for all $|t|<2L$, hence it holds for all $t$. {\hfill $\Box$ \vspace{2mm}}

Now let  $M$ always be a closed surface such that all geodesics are closed 
(unless otherwise stated). \\ The 
set $S:= \{  v\in SM \quad | \quad \gamma_v \quad \mbox{is simple closed } \}$
is a closed subset of $SM$ since the dimension of $M$ is $2$. Let $S^c$ denotes the complement of $S$ in $SM$. \\ 

The next lemma will be useful in the next section:
\begin{lem}
If $v_i\in S^c \to v\in S$ then $\liminf_{i\to \infty} P(v_i) \ge 2P(v)$.
\end{lem}
Proof: Assume given $0 \le t_{1,i}< t_{2,i} < P(v_i)$ such that $\gamma_{v_i}(t_{1,i}) = \gamma_{v_i}(t_{2,i})$. We have $|t_{1,i} - t_{2,i}| \ge inj(M)$, and 
$P(v_i)-|t_{1,i} - t_{2,i}| \ge inj(M)$. \\ Let $L=\liminf_{t \to \infty} P(v_i)$
and choose a subsequence (still denoted by $v_i$) such that $P(v_i) \to L$
and $ t_{1,i}\to s_1 $ and $ t_{2,i}\to s_2 $. Note that since $P$ is a discrete function
we have $P(v_i)=L$ for large $i$. Since $\gamma_{v_i}\to \gamma_v$
we have $\gamma_v(s_1) = \gamma_v(s_2)$. The assumption that $\gamma_v$ is simple closed
implies $|s_1-s_2| = nP(v)$ where $n\in \mathbb{N}$. Furthermore $\gamma_v(0)= \gamma_v(L)$
and therefore $L= mP(v)$ where $m\in \mathbb{N}$. This implies:
$$L= \lim P(v_i) \ge |s_1-s_2|+ inj(M) \ge P(v) +inj(M) $$
and therefore $L \ge 2P(v)$.{\hfill $\Box$ \vspace{2mm}}

The existence of simple closed geodesics on compact Riemannian manifolds $(M,g)$ of dimension 2 will play a crucial role for proving our main result. The following theorem garantees
us the existence of a sufficient large number of simple closed geodesics.

\begin{satz}[Ballmann]
Every compact Riemannian manifold $(M,g)$ of dimension 2  has at least three simple closed geodesics.
\end{satz}
Proof: See [B]. {\hfill $\Box$ \vspace{2mm}}

\begin{lem}
If $\gamma$ is a  closed geodesic on a surface $(M,g)$ and $(M,g)$ is a P-manifold
then every geodesic (except $\gamma$) intersects $\gamma$. 
\end{lem}

Proof: Take $x\notin \gamma$ and set $O=\{w\in S_xM \quad | \quad \gamma_w \quad\mbox{intersects}\quad \gamma    \}$ and $W=SM|\gamma-T\gamma$ (the set of all
 unit tangent vectors based on $\gamma$ which are not element of $T\gamma$). It is clear that $O$ is open, but it is also closed. Indeed,
choose a sequence $v_i\in O \to v$ with  $0 < t_i \le \max P $ 
and $\Phi_{t_i}(v_i)\in W$. Since $t_i$ is bounded, there exists a converging subsequence (still denoted by $t_i$). If $\Phi_{t_i}(v_i)\to \Phi_{s}(v) \notin W$, then $\Phi_{s}(v) \in T\gamma$ hence $x \in \gamma$. {\hfill $\Box$ \vspace{2mm}}

In the next section the Euclidean isometries of $\mathbb{R}^3$ are used to characterize the geodesic flows. The following theorem is important
\begin{satz}
If $F$ is a periodic  homeomorphism (i.e. $F^n=Id$ for some $n > 0$) of $S^2$ then $F$ is topologically  conjugate to 
the restriction of an Euclidean isometry of the ambient 3-space.
\end{satz}
Proof: See  [C] or  [K].{\hfill $\Box$ \vspace{2mm}}

\begin{ko}
Let $F$ be a periodic homeomorphism of $S^2$. If $F$ is orientation preserving
and has at least 3 fixed points then $F=Id$. If $F$ is orientation reversing and has at least one fixed point, then the fixed point set is an embedded circle.
\end{ko}
Proof: Theorem [1.9] implies that $F$ is conjugate to an Euclidean isometry. If $F$ is 
orientation preserving then $F$ is a rotation whose normal form is given by 
$$A = \left( \begin{array}{ccc}  \cos (\phi )  &  
\sin (\phi ) & 0 \\  -\sin (\phi ) & \cos (\phi ) & 0 \\ 0 & 0 & 1 \end{array} \right) .  
$$  
If this map has 3 fixed points then  the map is the identity, hence $F=Id$.\\ If $F$ 
is orientation reversing then $F$ is conjugate to a reflection whose normal form is given by 
$$A = \left( \begin{array}{ccc}  \cos (\phi )  &  
\sin (\phi ) & 0 \\  -\sin (\phi ) & \cos (\phi ) & 0 \\ 0 & 0 & -1 \end{array} \right) . 
$$ 
If this map has one fixed point then the normal form is given by
$$A = \left( \begin{array}{ccc}  -1  &  0 & 0 \\  0 & 1 & 0 \\ 0 & 0 & 1 \end{array} \right) .$$
The fixed point set is an embedded circle and therefore the fixed point set of $F$ is an embedded circle as well. {\hfill $\Box$ \vspace{2mm}}

\section[]{$\mathbb{R}P^2$ as a P-manifold}

In this section we prove our  main result, so 
let us assume that $M=\mathbb{R}P^2$ and $(M,g)$ is a P-manifold. In this case the
geodesic flow $\Phi$ is periodic (see theorem [1.4]) and equicontinuous (see lemma [1.5]).
A well-known technic to study a geodesic flow on a surface is to study the geodesic return map. 
Let $A= S^1\times (0,1) $ denote the open annulus. Let $\gamma$ be a simple closed geodesic. Denote by $W=SM|\gamma-T\gamma$ (the set of all
 unit tangent vectors based on $\gamma$ which are not element of $T\gamma$).
We want to prove by using the geodesic return map that $(M ,g)$ has infinitely many distinct simple closed geodesics. It is possible to define the geodesic return map on $M$, but we must separate two cases.  \\ Let us assume that $\gamma$ is a simple closed geodesic  that preserves the orientation.
Note that in this case $W$ is homeomorphic to an union of two  open annulus $A_0$ and $A_1$. We can identify $v \in A_0$ with $(x,\theta )$ where $\pi(v)=x$ and 
$\theta \in (0,1)$ is 
the angle of $v$ and $\dot{\gamma}(t)$ divided by $\pi$.  \\ Since every orbit is periodic we define  for our flow $\Phi$ the map  $F:A_0\to A_0 $  by 
$$ F(x,\theta )=(x_0, \theta_0) $$
where $x_0=\pi (\Phi_{t_0}(x,\theta))$   is the next intersection point of $\{ \pi (\Phi_t(x,\theta))|t>0 \} $ with 
$\gamma$   such that $\Phi_{t_0}(x,\theta)=(x_0, \theta_0)\in A_0$.
We just simple write $F:A\to A$. \\ If  $\gamma$ reverses the orientation, then
we find a neighbourhood of $\gamma$ looking like a M\"obius strip.
 Denote $W=SM|\gamma-T\gamma$ (the set of all
 unit tangent vectors based on $\gamma$ which are not element of $T\gamma$).
The set $W$ is homeomorphic to an open annulus, since $W$ is connected and orientable. 
Since every orbit is periodic (recurrent) we define  for our flow $\Phi$ the map  
$F:W\to W $  by 
$$ F(v)=v_0 $$
where $\Phi_{t_0}(v)=v_0 \in W$ and $x_0=\pi (\Phi_{t_0}(v))$   is the next intersection point of $\{ \pi (\Phi_t(v))|t>0 \} $. Since $W$ is homeomorphic to $A$, we  simple write
$F:A\to A$. 
\begin{lem}
Let $S^2$ be the two-point-compactification of $A$, i.e. $S^1\times \{0\} \sim -\infty$
and $S^1\times \{1\} \sim \infty$. Then $F$ can be extended to an homeomorphism $F:S^2\to S^2$
such that $\{ -\infty  , \infty \} $ are fixed points. 
\end{lem}
Proof: It is clear that $F:A\to A$ is an homeomorphism.  
If for $v=(x,\theta)$ we have $\theta$ near zero, the geodesic $\Phi_t(v)$
stays near $\dot{ \gamma }$ (equicontinuity) hence the second coordinate of  $F(x,\theta)$ is near zero, thus  $\{\infty\}$ and $\{- \infty\}$ 
are fixed points and moreover $F$ is continuous on $S^2$.{\hfill $\Box$ \vspace{2mm}}

We call  in this paper the  extension of $F$ to $S^2$ the geodesic return map
and denote it with $F:S^2\to S^2$. \\ Note that $S$ (see chapter 2) induces a set $$S_{A}:=\{v\in A \quad | \quad \gamma_{v} \quad \mbox{is simple closed} \} \subset A $$ in $A \subset S^2$, due to lemma [1.8] and theorem [1.7]. $F$ is a pointwise periodic homeomorphism on $S^2$ and therefore $F$ is periodic, i.e. $F^n=Id$ for some $n> 0$, since the following theorem holds:
\begin{satz}[Montgomery]
Let $F:X \to X$ be a pointwise periodic homeomorphism on a compact manifold $X$
then $F$ is periodic.
\end{satz}
Proof: See [M]. {\hfill $\Box$ \vspace{2mm}}

\begin{lem}
Given a simple closed geodesic $\gamma$ with length $L$ and the geodesic return map  $F:S^2\to S^2$ 
constructed from this geodesic. \\ Set $Per(x):=\inf\{n>0 \quad | \quad F^n(x)=x \} $ for $x\in A $. Then the following holds: \begin{enumerate}[i)]
\item If $\gamma$ preserves the orientation then 
$$ 2Per(v) = \sharp \{ \gamma_{v}([0, P(v))) \cap  \gamma([0, L))\}. $$
\item If $\gamma$ reverses the orientation then 
$$ Per(v) = \sharp \{ \gamma_{v}([0, P(v))) \cap  \gamma([0, L))\}. $$
\end{enumerate}
\end{lem}

Proof: ii) is clear by construction and therefore we only prove i).\\ We take the model of $M= \mathbb{R}P^2$
as a closed disk where we identify the boundary points. Set the boundary curve to be $\beta$.
Let $V$ be a vectorfield along $\gamma$ with $g(V,V)=1$ and 
$g(V(t), \dot{\gamma}(t))= 0$. Therefore $V(t)=V(t+P_0)$, where $P_0$ is the period of $\gamma$. \\ We define  what lies above 
and below near $\gamma$. Given a small neighbourhood $U$ of $\gamma$, so small that $U$ is diffeomorphic to $S^1\times (0,1)$. $U-\gamma$ is an union of two open annulus $U_0$ and $U_1$. Let $U_0$ denote the set where the vectors of $A$ are pointing in. We say that a point in $U_0$ lies above and a point in $U_1$ lies below.
Fix a closed geodesic $c$ that cuts $\gamma$ transversally. We show that $c$ intersects $\gamma$
 mutually from above and below.
If this holds  for a small pertubation of $\gamma$ in the space of the $C^{\infty}$ simple closed curves then it holds for $\gamma$. Thus we can suppose that $\gamma$ intersects $\beta$ and $c$  transversally. It is easy to see that $\gamma$ intersects $\beta$ for an even times, otherwise it would reverses orientation.  
Thus we can lift $\gamma$ to a simple closed curve on $S^2$. We lift the vectorfield $V$ up and conclude  that $c$ intersects $\gamma$ as expected (regard a lift of $c$ and use the Jordan curve theorem).{\hfill $\Box$ \vspace{2mm}}

\begin{lem}  
Given the geodesic return map $F:S^2\to S^2$ induced from $\gamma$. \\Set  $m_0=\min\{ Per(x) | x\in A  \} $. Then 
 the following holds: 
\begin{enumerate}[i)]
\item $F^{m_0}$ has three fixed points
\item $Per(A)=\{  m_0, 2m_0 \}$
\item If $F^{m_0}$ is orientation preserving then $F^{m_0}=Id$. 
      Furthermore all geodesics are simple closed.
\item If $F^{m_0}$ is orientation reversing the fixed point set is an embedded circle.
\item If for some $v\in S_A$ we have $Per(v)=2m_0$ 
      then $v$ is an interior point of $S_A$, i.e. $S_A$ contains an open set.
\end{enumerate}

\end{lem}
Proof: i) is clear, since $\{ \pm \infty \} $ are fixed points.\\ ii)
If $m_0 < Per(y)=p < 2m_0$ then $F^{2m_0-p}(y)=y$ and \\ $0 < 2m_0-p< m_0$. \\ iii)
Since $F$   has three fixed points,
we conclude from corallary [1.10] that $F=Id$. We show that all geodesics are simple closed.
We first prove the case when $\gamma$ preserves orientation.\\ Lemma [2.3, i)] shows that all geodesics (except $\gamma$) are intersecting $\gamma$ for  $2m_0$ times, i.e. $m_0\equiv  Per$. We show that $S_A$ is an open set.
Assume there is a sequence $v_i\to v\in S_A$ such that $\gamma_{v_i}$ has selfintersections.
It follows from lemma [1.6] that we can suppose $P(v_i) \ge  nP(v)$ where $n>1$. If $i$ grows it follows from transversality that $\gamma_{v_i}$ intersects 
$\gamma$ more than $2m_0$ times. Take a look at the picture.

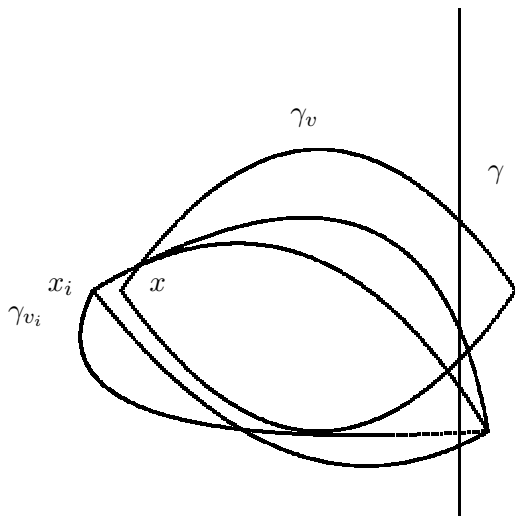
\setlength{\unitlength}{0,75cm}
\begin{picture}(12,10)
\linethickness{1pt}
\qbezier[200](2.5,5)(6,0.5)(9.5,2.5)
\qbezier[200](2.5,5)(1,2)(9.5,2.5)
\qbezier[200](9.5,2.5)(8.5,8.5)(2.5,5)
\qbezier[200](9.5,2.5)(6.5,7.5)(2.5,5)
\qbezier[200](3,5)(6.5,0)(10,5)
\qbezier[200](10,5)(6.5,10)(3,5)

\put(3.5,5){$x$}
\put(1.7,5){$x_i$}
\put(9.5,7){$\gamma$}
\put(6,8){$\gamma_{v}$}
\put(1,4.5){$\gamma_{v_i}$}

\thicklines
\put(9,10){\line(0,-1){9}}

\end{picture}

Regard a local coordinate system around $\gamma$. Choose a sequence \\ $0\le t_{i,1}< t_{2,i} < P(v_i)$  such that $\gamma_{v_i}(t_{1,i}) = \gamma_{v_i}(t_{2,i})$. Set $x_i=\gamma_{v_i}(t_{1,i})$. Note that $\dot{\gamma}_{v_i}(t_{1,i}) $ and 
$\dot{\gamma}_{v_i}(t_{2,i})$ converge  to  
tangent vectors  of $\gamma_v$. \\ W.l.o.g. $x_i \to x \notin \gamma$. Thus if $i$ grows,  we have for a small $\delta > 0$  that
$t_{1,i}+P(v_0) -\delta < t_{2,i}+ \delta$ . Moreover, we can assume that
$\gamma_{v_i}([t_{1,i}+\delta ,t_{1,i}+ P(v_0)-\delta ]  ) $ and 
$ \gamma_{v_i}([t_{2,i}+\delta ,t_{2,i}+ P(v_0)-\delta ]  ) $ intersect $\gamma$ for $2m_0$ times, thus $\gamma_{v_i}([ 0, P(v_i)] )$ intersects $\gamma$ for at least  $4m_0$ times. This gives a contradiction.
Hence $S_A$ is open and nonempty, but it
 is closed in $A$.\\ If $\gamma$ is  orientation reversing the proof is analogue. One must apply lemma [2.3, ii)] instead of lemma [2.3, i)]. \\ iv) If $F^{m_0}$ is orientation reversing one can deduce as in the proof of lemma [2.3, iii)] from corollary [1.10]  that its fixed point set is an embedded circle.\\ v) We only prove the case where $\gamma$ is orientation preserving, since the other case is analog. All geodesics (except $\gamma$) are intersecting $\gamma$  at most $4m_0$ times! Assume there is  a sequence $v_i \in S_A^c \to v \in S_A$ where $Per(v)=2m_0$. Since $P(v_i)\to nP(v)$ where $n> 1$ we conclude as in the proof of lemma [2.4, iii)] that $\gamma_{v_i}$ 
intersects $\gamma $ more than $4m_0$ times 
for $i$ sufficiently large.{\hfill $\Box$ \vspace{2mm}}

\begin{hl}
Given   a simple closed geodesic $\gamma$ and the associated geodesic return map
$F:S^2\to S^2$.
Then  $|S_A|= \infty$, hence there are infinitely many 
distinct simple closed geodesics.
\end{hl}

Proof: Let us set $G:= F^{m_0}$. W.l.o.g. $G $ is orientation reversing, otherwise we apply lemma [2.4, iii)].  By corollary [1.10]  we have  $F^{2m_0}= G^2 =Id$, since $G^2$ has at least three fixed points
and is orientation preserving.  \\ We first regard the case where $\gamma$ preserves the orientation (i.e. we can find a  neighbourhood of 
$\gamma$ homeomorphic to a cylinder). Assume $|S_A|$  is finite thus $Per(S_A)=m_0$ by lemma [2.4, v)] and  $G:= F^{m_0}$ has at least $2+2$ fixed points. Indeed, we have 
$$Fix(G)= \{v\in A \quad | \quad Per(v)=m_0 \} \cup \{\pm \infty \}, $$ and the remaining two simple closed geodesics given by theorem [1.7] and $\{\pm \infty \}$  gives us at least $4$ fixed points  by the previous results.  $G$  is orientation reversing 
and therefore by corollary [1.10]  its fixed point set is an embedded circle. We conclude from the method of lemma [2.4, iii)] that $S_A$ 
contains an arc connecting $\infty $ and $-\infty $. Indeed,  $Fix(G)-\{ \infty , -\infty \} $ is the union of two arcs that contains at least two fixed points (that induce simple closed geodesics). There is an arc $\nu \subset Fix(G)-\{ \infty , -\infty \}$ that contains at least one point from $S_A$
and is connecting $\infty $ and $-\infty $.
 We show that  $S_A \cap \nu $ is open in the space $\nu$! \\ Take $v_i\in S_A^c\cap \nu \to v\in S_A \cap \nu$, therefore we have 
$Per(v)=m_0=Per(v_i)$. Again $P(v_i)\to nP(v)$ where $n> 1$ and 
we conclude as in lemma [2.4, iii)] that $Per(v_i)> m_0$ if $i$ grows, which 
gives a contradiction.
Therefore $S_A\cap \nu$ is an open set in $\nu$. \\ If $\gamma$ reverses the orientation
then the proof is analoug. 
{\hfill $\Box$ \vspace{2mm}}

Now we are able to prove our  theorem:
\begin{satz}
Let $g$ be a  metric on $\mathbb{R}P^2$. Then $(\mathbb{R}P^2 , g)$
has constant curvature iff all geodesics are closed. 
\end{satz}

Proof: We  show that $(\mathbb{R}P^2, g)$ is a SC-manifold, since 
this implies that $(\mathbb{R}P^2, g)$  has constant curvature (see 10.10.3 theorem 258 in [Be]).
Let us denote $M= \mathbb{R}P^2 $ and $SM$ its unit tangent bundle. We know that the geodesic flow $\Phi$ on $SM$ 
is a smooth semi-free (sometimes called locally-free) $S^1$-action. Suppose that 
$u \in SM$ is a point such that the orbit of $\Phi$ does not have maximal period, i.e. the isotropy group  $S^1_u$ at $u$
 is a non-trivial cyclic group generated by say 
$f_k = e^{\frac{2\pi i}{k}} \in S^1$. Note that $SM$ as a contact manifold is orientable and the finite order diffeomorphism $f_k: SM \to SM$ is orientation preserving.  Thus its fixed point set is one dimensional and hence exactly the orbit
through $u$.  Any orbit of this type (that has no maximal period) is isolated and  therefore there are at most
finitely many such exceptional orbits. Let us prove this in detail: \\ Given   $SM$ and a smooth circle action of $SM$, then we can find a smooth Riemannian metric $g_0$ invariant under the circle action by
averaging of a smooth  Riemannian metric on $SM$. Thus any  element of our action  becomes
an isometry of $(SM,g_0)$.  Given now $f_k$ and $u\in Fix(f_k)$.
Note that the orbit of $u$ under the  geodesic  flow $\Phi$ lies in $Fix(f_k)$. Suppose that  $Fix(f_k)$ is not the union of
finitely many such exceptional orbits, therefore
for some point $x\in Fix(f_k)$ we find for $Df(x)$ 2 independent eigenvectors with eigenvalue $1$. Take a look at the picture that shows how $Fix(f_k)$ looks likewise if it is not a manifold of dimension 1.

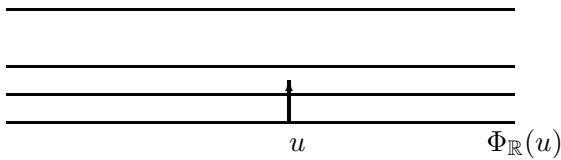
\setlength{\unitlength}{0,75cm}
\begin{picture}(12,4)
\linethickness{1pt}

\put(5,0.5){$u$}
\put(8.5,0.5){$\Phi_{\mathbb{R}}(u)$}
\put(5,1){\vector(0,1){0,75}}
\thicklines
\put(0,1){\line(1,0){9}}
\put(0,1.5){\line(1,0){9}}
\put(0,2){\line(1,0){9}}
\put(0,3){\line(1,0){9}}

\end{picture}

Note that $(Df_k(x))^k=Id$ (thus $Df_k(x)$ has no nilpotent part) and therefore we conclude $f_k=Id$, since $f_k: SM \to SM$ is orientation preserving and an isometry of $(SM, g_0)$ such that $Df_k(x)=Id$. This implies that $Fix(f_k)$  is near $u$  a segment of the  unique orbit through $u$  and therefore the finite union of closed orbits. The same arguments can be applied to $f_k^l$ 
where $l\in\{2, \cdots, k-1 \}$. \\ Note that $S$ is the set of all $v\in SM$ that induces simple closed geodesics. Let us denote by $E$  the dense, open and 
connected set in $SM$ containing only the no-exceptional orbits  ($SM$ has dimension three). Note that $S^c \cap E$, as a subset of $E$, is open and there must exist an orbit  $v \in E$ that is simple closed by proposition [2.5], since $E^c$ is the union of finitely many orbits. The set $S \cap E$ is open in $E$.
Otherwise we approximate $v\in S \cap E$ with $v_i \in S^c\cap E$ and conclude
that $v$ has selfintersections by lemma  [1.6]. Thus $E = E \cap S \subset S$ and since $S$ is closed in $SM$, 
we conclude $S=SM$. This argument was also used in [G].\\  Finally  we show that the period  for all geodesics are the same.
Take a geodesic $\gamma$ that has maximal period and construct the map 
$F:A\to A$.  All exceptional points $E^c$ induce a set $E^c_A= E^c\cap A$, since lemma [1.8] holds. We prove only the case where $\gamma$ preserves orientation, since the
other case is analogue. Note $|E^c_A|< \infty $, since $E^c$ is the union of finitely many  orbits. Given any $w\in E^c_A$ that has no maximal period. Let $P_0$ denote the maximal period. Take a sequence $w_i \in E_A=E\cap A \to w $ such that $Per(w_i)=Per(w)$ 
and $bP_0=bP(w_i)=aP(w)$, where $a,b$ are  integers. This is possible since $Per$ is constant on $A$ or $Fix( F^{m_0}) $ is an embedded circle and the set of those points $v$ such that $Per(v)\not= m_0$ is open. The closed curves $\gamma_{w_i}:[0, bP_0] \to M $ converge to 
$\gamma_{w}:[0, aP(w)] \to M $, thus, if $i$ grows, the intersection number of these closed curves   with $\gamma$ is the same as that of $\gamma_w$ with $\gamma$. Therefore $2bPer(w_i)= 2aPer(w)$ 
and this shows $a=b$ and so $P(w)=P_0$.
{\hfill $\Box$ \vspace{2mm}}

As a corollary we conclude:
\begin{satz}
All two dimensional P-manifolds are SC-manifolds.
\end{satz}
Proof: The only P-manifolds in dimension 2 are $S^2$ and $\mathbb{R}P^2$ (see theorem 7.37 in [Bs]). The first case follows from 
[G] and the second case follows from theorem [2.6].{\hfill $\Box$ \vspace{2mm}}

\section[]{Zoll surfaces }
In this small section we want to point out the remarkable difference between $P$-metrics on 
$S^2$ and $\mathbb{R}P^2$. We state the results of Zoll and Darboux that easily shows the existence
of many different $P$-metrics on $S^2$. The following classes of Riemannian manifolds $(S^2,g)$ are important in the theory of 2-dimensional $P$-manifolds:

\begin{defin}
A metric $g$ on $S^2$ is called a Zoll metric if g is continuous and $(S^2,g)$
is a $C_{2\pi}$-surface (i.e. all geodesics are periodic and have the  period $2\pi$).
In this case $(S^2,g)$ is called a Zoll surface. 
\end{defin}

\begin{defin}
A metric $g$ on $S^2$ is called a  metric of revolution if $(S^2,g)$ has $S^1$ as an effective isometry subgroup. The surface $(S^2,g)$ is called a surface of revolution.
\end{defin}

 If $(S^2,g)$ is a surface of revolution then the action $S^1$ has exactly two fixed points
, say $N$ and $S$ (the so-called North and South poles). \\ Since $ U:=S^2- \{S,N \} $ is diffeomorphic
to the punctured disc $$\{ z\in \mathbb{C} \quad | \quad |z|<1, z \not= 0 \},$$ we can find
 a coordinate system such that $g$ can be written on $U$ as 
$$g= du^2  + a^2(u) \theta^2$$
such that $a^2$ is bounded by $1$. From this we can find a coordinate system
$[U, (r,\theta )]$ such that  $g$ can be written on $U$ as 
$$g= f(cosr)^2dr^2  + sin^2rd\theta^2,$$
where $f$ is a function from $(-1,1)$ to $\mathbb{R}^{+}$  (compare proposition 4.10 in [Bs]). The following theorem holds:

\begin{satz}[Darboux]
Let $(S^2,g)$ be a surface of revolution and $g$ be written in the parametrization $[U, (r, \theta )]$ as 
$$g= f(cosr)^2dr^2  + sin^2rd\theta^2. $$
A necessary and sufficient condition in order that all geodesics are closed is that
for every $t$ one has
$$ \int_t^{\pi - t} \frac{f(\cos r) \sin t  }{\sin r \sqrt{\sin^2r - \sin^2t}}\, dr
= \frac{p}{q}\pi $$
where $p$ and $q$ are integers.
\end{satz}
Proof: See [D] or theorem 4.11 in [Bs]. {\hfill $\Box$ \vspace{2mm}}

If we assume that $(S^2,g)$ is a Zoll surface of revolution then the following theorem holds:

\begin{satz}
A metric $g$ is a Zoll metric of revolution on $S^2$ iff $g$ can be written in the parametrization $[U, (r, \theta )]$ as 
$$g= (1+h(\cos r))^2\, dr^2+\sin^2r\theta^2\, d\theta^2 $$
where $h$ is an odd function from $[-1,1]$ to $(-1,1)$, whose value at $1$ is $0$.
Furthermore, $(S^2,g)$ is a $SC_{2\pi}$-manifold (i.e. a $C_{2\pi}$-manifold whose geodesics are all simple closed). Moreover, $g$ is $C^k$ at $N$ (resp. $S$) iff $h$ extends to a $C^{k-1}$ function in a neighbourhood of $1$ (resp. $-1$).
\end{satz}
Proof: See theorem 4.13 and corollary 4.16 in [Bs].{\hfill $\Box$ \vspace{2mm}}

It is remarkable that the choice of the function $h$ does not involve any condition on the derivative. Since the sectional curvature can be expressed in the form 
$$\sigma (r) = \frac{1}{(1+h(\cos r))^3} (1 + h(\cos r)-\cos r\cdot h'(\cos r)) ,$$
we see that there is a large class of $C^{\infty}$ Zoll metrics whose sectional curvature is not constant. One can choose for example \\ $h(\cos r)= \cos r\cdot \sin(2k+1)r$ or $h(\cos r)=\cos r \cdot \frac{sin^2r}{2}$. In paricular, one can prove that there is a large
set of different $P$-metrics.

\section[]{Acknowledgements}

The author thanks Gerhard Knieper, who motivated the author to research on the field
of P-manifolds and Karsten Grove, Uwe Abresch and Henrik Koehler for their discussions
and corrections.

Christian Pries \\ Fakult$\ddot{a}$t f$\ddot{u}$r Mathematik\\ Ruhr-Universit$\ddot{a}$t 
Bochum\\ Universit$\ddot{a}$tsstr. 150 \\ 44780 
Bochum \\ Germany\\ Christian.Pries@rub.de

\end{document}